\documentclass[12pt]{article}
\usepackage{amsmath}
\usepackage{latexsym}
\usepackage{amssymb}
\newtheorem{thm}{Theorem}[section]
\newtheorem{la}[thm]{Lemma}
\newtheorem{Defn}[thm]{Definition}
\newtheorem{Remark}[thm]{Remark}
\newtheorem{Note}[thm]{Note}
\newtheorem{prop}[thm]{Proposition}

\newtheorem{Example}[thm]{Example}
\newtheorem{Examples}[thm]{Examples}
\newtheorem{Problems}[thm]{Problems}

\newtheorem{Problem}[thm]{Problem}
\newtheorem{Number}[thm]{\!\!}

\newenvironment{example}{\begin{Example}\rm}{\end{Example}}

\newenvironment{rem}{\begin{Remark}\rm}{\end{Remark}}
\newenvironment{numba}{\begin{Number}\rm}{\end{Number}}
\newenvironment{proof}{{\noindent\bf Proof.}}%
                  {\nopagebreak\hspace*{\fill}$\Box$\medskip\medskip\par}   
\newcommand{\Punkt}{\nopagebreak\hspace*{\fill}$\Box$}

\newcommand{\wb}{\overline}

\newcommand{\wt}{\widetilde}

\newcommand{\impl}{\Rightarrow}

\newcommand{\mto}{\mapsto}

\newcommand{\isom}{\cong}

\newcommand{\N}{{\mathbb N}}
\newcommand{\R}{{\mathbb R}}
\newcommand{\bD}{{\mathbb D}}

\newcommand{\K}{{\mathbb K}}

\newcommand{\Z}{{\mathbb Z}}
\newcommand{\C}{{\mathbb C}}

\newcommand{\cO}{{\mathcal O}}

\newcommand{\dl}{{\displaystyle \lim_{\longrightarrow}}}

\newcommand{\sub}{\subseteq}

\newcommand{\cH}{{\mathcal H}}

\DeclareMathOperator{\Spann}{span}

\begin{document}
\begin{center}
{\Large\bf
Instructive examples of smooth,\\[.9mm]
complex differentiable and complex analytic\\[2.7mm]
mappings into locally convex spaces}\\[5mm]
{\bf Helge Gl\"{o}ckner}\vspace{2mm}
\end{center}
\begin{abstract}
\hspace*{-7.2 mm}
For each $k\in \N$,
we describe a mapping $f_k\colon \C\to E_k$
into a suitable non-complete complex
locally convex space~$E_k$
such that $f_k$ is
$k$ times continuously complex differentiable
(i.e., a $C^k_\C$-map)
but not $C^{k+1}_\C$
and hence not complex analytic.
We also describe a complex analytic
map from~$\ell^1$
to a suitable complete complex locally
convex space~$E$ which is unbounded
on each non-empty open subset of~$\ell^1$.
Finally, we present
a smooth map $\R\to E$ into a
non-complete locally convex space which is
not real analytic although it is given locally
by its Taylor series around each point.
\end{abstract}
{\footnotesize {\em Classification}:
46G20 (primary), 
26E05, 
26E15, 
26E20, 
46T25\\[1mm] 
{\em Key words}:
complex differentiability,
smoothness,
analyticity,
analytic map,
holomorphic map,
infinite-dimensional holomorphy,
infinite-dimensional calculus,
local boundedness,
completeness, free locally convex space}\vspace{3mm}
\begin{center}
{\Large\bf Introduction}
\end{center}
It can be advantageous to
perform infinite-dimensional
differential calculus in
general locally convex spaces,
without completeness
conditions.
First of all,
the theory becomes clearer and more transparent
if completeness conditions
are stated explicitly as hypotheses
for those results which really depend
on them, but omitted otherwise.
Secondly,
it simplifies
practical applications if
completeness properties
only need to be checked
when they are really needed.
Therefore, various authors
have defined and discussed $C^k$-maps
(and analytic maps)
between locally convex spaces
without completeness
hypotheses
(see
\cite{RES}, \cite[Chapter~1]{GaN},
\cite{Gro} and \cite{Kel};
cf.\ \cite{BaS}).\\[2.5mm]
In this article,
we compile examples
which illustrate
the differences between various
differentiability and
analyticity properties of
vector-valued functions,
in particular differences
which depend on non-completeness
of the range space.
Primarily, we consider continuous
mappings $f\colon U\to F$,
where $U\sub \C$
is open and~$F$
a complex locally convex space.
Let us call such a map~$C^1_\C$
if the complex
derivative $f^{(1)}(z)=f'(z)=\frac{df}{dz}(z)$
exists for each $z\in U$, and $f'\colon U\to F$
is continuous.
As usual, we say that~$f$ is $C^k_\C$
if it has\linebreak
continuous
complex derivatives $f^{(j)}\colon U\to F$
for all $j\in \N_0$ such that $j\leq k$
(where $f^{(j)}:=(f^{(j-1)})'$).
Finally, call~$f$ \emph{complex analytic}
if it is of the form
$f(z)=\sum_{n=0}^\infty (z-z_0)^na_n$
close to each given point $z_0\in U$,
for suitable elements $a_n\in F$.
The weakest relevant
completeness property of~$F$
is \emph{Mackey completeness},
which (among many others) can be defined by
each of the following
equivalent conditions:\,\footnote{See \cite[Theorem~2.14]{KaM},
also \cite[p.\,119]{Kri}.}
\begin{itemize}
\item[M1]
The Riemann
integral $\int_0^1 \gamma(t)\,dt$
exists in~$F$ for each smooth curve
$\gamma\colon \R\to F$.
\item[M2]
$\sum_{n=1}^\infty t_n x_n$
converges in~$F$ for each bounded
sequence $(x_n)_{n\in \N}$
in~$F$
and each sequence~$(t_n)_{n\in \N}$
of scalars such that $\sum_{n=1}^\infty |t_n|<\infty$.
\end{itemize}
If $F$ is Mackey complete,
the following properties
are known to be equivalent:\footnote{See Theorems~2.1,
2.2 and~5.5 in \cite[Chapter~II]{Gro},
or simply replace sequential completeness
with Mackey completeness
in~\cite[Theorem~3.1]{BaS}
and its proof.}
\begin{itemize}
\item[(a)]
$f\colon \C\supseteq U\to F$ is complex analytic;
\item[(b)]
$f$ is $C^\infty_\C$;
\item[(c)]
$f$ is $C^1_\C$;
\item[(d)]
$f$ is weakly analytic,
i.e.\ $\lambda\circ f\colon U\to\C$
is complex analytic for each
continuous linear functional $\lambda\colon F\to\C$;
\item[(e)]
$\int_{\partial\Delta} f(\zeta)\,d\zeta=0$
for each triangle $\Delta\sub U$;
\item[(f)]
$f(z)=\frac{1}{2\pi i}\int_{|\zeta-z_0|=r}\frac{f(\zeta)}{\zeta-z}\,d\zeta$
for each $z_0\in U$,
$r>0$ such that
$z_0+r\bD\sub U$,
and each~$z$
in the interior of $z_0+r\bD$
(where $\bD:=\{z\in \C\colon |z|\leq 1\})$.
\end{itemize}
Here (a) and (b) remain equivalent
if~$F$ fails to be Mackey complete
(see \cite[Chapter~II, Theorem~2.1]{Gro}
or \cite[Propositions~7.4 and~7.7]{BGN})
and also~(d), (e) and~(f) remain
equivalent (because
$f$ may be considered as a map
into the completion of~$F$
for this purpose; see also \cite[Chapter~II,
Theorem~2.2]{Gro}).
However, $C^1_\C$-maps need not be
$C^2_\C$ (and hence
need not be complex analytic) in this case,
as an example in
Herv\'{e}'s book~\cite[p.\,60]{Her} shows
(for which only a partial proof is provided there).\\[2.5mm]
Our first goal is to give
examples which distinguish
between the properties (a)--(f) in a
more refined way.
Thus, we describe functions
satisfying~(d)--(f) but which are not~$C^1_\C$,
and also $C^k_\C$-maps which are not $C^{k+1}_\C$,
for each $k\in \N$ (Theorem~\ref{mainres}).
In particular, the latter functions are
$C^1_\C$ but not complex analytic,
like Herv\'{e}'s example.\\[2.5mm]
We mention that similar functions
have been recorded
in unpublished parts of the thesis~\cite{Gro}
(Chapter~II, Example~2.3),
but also there a crucial step
of the proof is left to the reader.
Furthermore, our discussion is based on a different
argument: At its heart is the linear
independence of the functions
$\N\to \C$, $n\mto n^kz^n$
in the space~$\C^\N$
of complex sequences, for $k\in \Z$
and $z\in \C^\times:=\C\setminus\{0\}$
(see, e.g.\ \cite[Lemma~1]{LaP}).
As a byproduct,
this argument
entails refined results
concerning completeness properties
of free locally convex spaces
over subsets of~$\C$ (see Proposition~\ref{consfree}),
which go beyond
the known general facts
concerning free locally convex spaces (as in~\cite{Usp}).\\[2.5mm]
Examples
concerning real analyticity
are given as well,
as are examples\linebreak
concerning maps $f\colon U\to F$,
where $E$ and~$F$ are locally convex spaces
over $\K\in \{\R,\C\}$
and $U\sub E$ an open set.
Various differentiability and analyticity
properties
of such maps~$f$ are
regularly used in applications
of infinite-dimensional calculus,
notably in infinite-dimensional Lie theory:
\begin{itemize}
\item[(i)]
Given $k\in \N_0\cup\{\infty\}$,
$f$ is called a \emph{$C^k_\K$-map}
in the sense
of Keller's $C^k_c$-theory
if $f$ is continuous,
the iterated directional
(real or complex) derivatives
\[
d^jf(x,v_1,\ldots, v_j)\; :=\;
(D_{v_j}\cdots D_{v_1}f)(x)
\]
exist in~$F$ for all $j\in \N$ such that
$j\leq k$, $x\in U$ and
$v_1,\ldots, v_j\in E$,\linebreak
and the maps $d^jf\colon U\times E^j\to F$
so defined are continuous
(see, e.g.,  \cite{Kel},
\cite{RES} and \cite{GaN}).
As usual, $C^\infty_\R$-maps are also
called \emph{smooth}.
\item[(ii)]
If $\K=\C$, $f$ is called \emph{complex analytic}
if $f$ is continuous
and for each $x\in U$, there exists a $0$-neighbourhood
$Y\sub U-x$ and continuous,
complex homogeneous polynomials $p_n\colon E\to F$
of degree~$n$ such that
%
\begin{equation}\label{sumhomo}
f(x+y)\; =\; \sum_{n=0}^\infty p_n(y)
\end{equation}
for all $y\in Y$,
with pointwise convergence (see \cite{BaS}).
\item[(iii)]
If $\K=\R$,
following \cite{Mil},
\cite{RES} and~\cite{GaN},
the map
$f\colon U\to F$
is called \emph{real analytic}
if it extends to a complex
analytic $F_\C$-valued
map on an open neighbourhood
of~$U$ in the complexification $E_\C$.
\end{itemize}
It is known
that complex analytic maps coincide with $C^\infty_\C$-maps
(see \cite[Propositions~7.4 and~7.7]{BGN}
or \cite[Chapter~1]{GaN}),
and furthermore compositions
of composable $C^k_\K$-maps
(resp., $\K$-analytic maps) are $C^k_\K$
(resp., $\K$-analytic) for $\K\in \{\R,\C\}$
(see \cite[Proposition~4.5]{BGN}
or \cite[Chapter~1]{GaN};
cf.\ \cite{RES}). 
Occasionally,
authors use a different notion
of real analytic maps:
\begin{itemize}
\item[(iv)]
In~\cite{BaS},
$f$
is called real analytic if it
is continuous
and locally
of the form~(\ref{sumhomo})
with continuous, real homogeneous
polynomials~$p_n$.
\end{itemize}
This approach is not useful in general
because one is not
able to prove (without extra assumptions)
that compositions of such maps
are again of the same type.
Finally:
\begin{itemize}
\item[(v)]
Bourbaki~\cite{BVA}
defines analytic maps
from open subsets
of normed spaces
to quasi-complete
locally convex spaces
in a way (not recalled here)
which forces such maps
to be locally bounded
(cf.\ \cite[3.3.1\,(iv)]{BVA}).\footnote{Local
boundedness follows from
the definition of $\cH_R(E_1,\ldots, E_n;F)$
in \cite[3.1.1]{BVA}.}${}^{,}$\footnote{The booklet~\cite{BVA}
is still of relevance for Lie theory
as it provides the framework for
the influential
and frequently cited
volume~\cite{BLG}.
Also, it is a rare source of information
on analytic maps on open subsets
of infinite-dimensional
normed spaces over valued fields.}
\end{itemize}
It is known that
Bourbaki's notion
of complex analyticity coincides
with the one from~(ii)
for maps between Banach spaces
(cf.\ \cite[Proposition~5.1]{BaS}).
It is also known that real
analyticity as in~(iv)
coincides with the one from~(iii)
if~$E$ is a Fr\'{e}chet space
and~$F$ Mackey complete
(cf.\ \cite[Theorem~7.1]{BaS}
for the crucial case where~~$F$
is sequentially complete).
Nonetheless, the concepts differ
in general.
We illustrate the differences
with simple examples:\\[2.5mm]
In Proposition~\ref{propreal}
(announced in \cite[\S1.10]{GCX}),
we describe a smooth map
$f\colon \R\to F$
to some (non-Mackey complete)
real locally convex space~$F$
which is locally given by its Taylor
series around each point (and thus
real analytic in the inadequate sense of~(iv)),
but not real analytic in the sense
of~(iii).\\[2.5mm]
In Section~\ref{secnbd}, we provide an example of a map
$f\colon \ell^1\to \C^\N$ on the space $\ell^1$
of absolutely summable complex sequences
which is complex analytic in the sense
of~(ii) but unbounded on each non-empty
open subset of~$\ell^1$.
As a consequence, $f$ is not
complex analytic in Bourbaki's sense
(as in~(v)).
In particular, this means that the equivalence
of~(i) and~(ii) in \cite[3.3.1]{BVA}
is false
(which asserts
that complex analytic maps
in Bourbaki's sense
coincide with complex differentiable
maps).\\[2.5mm]
For complex analytic maps
in the sense of~(ii),
\emph{ample} boundedness
can be used as an
adequate substitute
for ordinary boundedness
(see \cite[Chapter~II, \S6, notably Theorem~6.1]{Gro};
cf.\ also \cite[Theorem~6.1\,(i)]{BaS}).\\[3mm]
{\bf General conventions.}
We write $\N=\{1,2,\ldots\}$
and $\N_0:=\N\cup\{0\}$.
If $(E,\|.\|)$ is a normed
space over $\K\in \{\R,\C\}$
(e.g., $(E,\|.\|)=(\K,|.|)$),
we write
$B^E_r(x):=\{y\in E\colon |y-x|<r\}$
and $\wb{B}^E_r(x):=\{y\in E\colon |y-x|\leq r\}$
for $x\in E$ and $r>0$.
Given a vector space~$E$ over a field~$\K$
and a subset $M\sub E$, we write $\Spann_\K(M)$
for the vector subspace of~$E$ spanned
by~$M$.
%
%
%
%
%
%
%
%
%
%
\section{Examples of {\boldmath$C^k_\C$}-maps which are not
{\boldmath$C^{k+1}_\C$}}\label{secmain}
For $k\in \Z$ and $z\in \C^\times$, we define
%
\begin{equation}\label{defhs}
h_{k,z}\colon \N\to \C\,,\quad h_{k,z}(n):=n^k z^n\,.
\end{equation}
Let $U\sub \C$ be a non-empty open subset
and $M \sub \C^\times$ be a superset of~$e^U$
(for example, $U=\C$, $M=\C^\times$).
For $k\in \N_0$, we let
$E_k\sub \C^\N$ be the vector subspace spanned
by the functions $h_{j,z}$ with $z\in M$ and
$j\in \N_0$ such that $j\leq k$.
We give $E_k$ the topology induced
by the direct product~$\C^\N$ and define
\[
f_k \colon U \to E_k\,,\quad f_k(z):=(e^{nz})_{n\in \N}\,.
\]
%
%
\begin{thm}\label{mainres}
For each $k\in \N$, the map
$f_k\colon U\to E_k$ is $C^k_\C$ but not
$C^{k+1}_\C$ $($and hence not
complex analytic$)$.
Furthermore,
$f_0\colon U \to E_0$ is
weakly analytic but not~$C^1_\C$
$($and hence not complex analytic$)$.
\end{thm}
The following fact
is crucial for our proof
of Theorem~\ref{mainres}.
It is known in the theory
of linear difference equations
with constant coefficients
(see, e.g., \cite{LaP}; cf.\ \cite{Nor}).
For the reader's convenience,
we give a self-contained proof.\footnote{An earlier
preprint version
gave a more complicated
proof, based on harmonic analysis.
The new proof only uses linear algebra.
It is a variant of an argument
communicated to the author by L.\,G. Lucht.
Note that $\C$ can be replaced
by any field of characteristic~$0$.}
%
\begin{la}\label{linind}
The functions $h_{k,z}$ $($for $k\in \Z$, $z\in \C^\times)$
are linearly
independent in $\C^\N$.
\end{la}
\begin{proof}
Step~1: To prove Lemma~\ref{linind},
we only need to show that
$(h_{k,z})_{k\in \N_0,z\in \C^\times}$ is a linearly
independent family of functions
in~$\C^\N$.
This follows from the fact
that the multiplication operator
$\C^\N\to \C^\N$, $f\mto h_{\ell,1}\cdot f$
with $(h_{\ell,1}\cdot f)(n)=n^{\ell}f(n)$
is a linear automorphism
and $h_{\ell,1}\cdot h_{k,z}=h_{k+\ell,z}$,
for all $\ell,k\in \Z$.\\[2.5mm]
Step~2: For each fixed
$z\in \C^\times$, the functions
$h_{k,z}$ ($k\in \N_0$)
are linearly independent.
In fact, if $m\in \N$,
$a_0,\ldots, a_m\in \C$ and
$\sum_{k=0}^m a_kh_{k,z}=0$,
then $(\sum_{k=0}^m a_kn^k)z^n=0$
for each $n\in \N$ and thus
$\sum_{k=0}^m a_kn^k=0$,
entailing that $a_0=a_1=\cdots=a_m=0$.\\[2.5mm]
Step~3: Now consider the shift operator
$S\colon \C^\N\to\C^\N$, $S(f)(n)=f(n+1)$.
Then $S$ is a linear endomorphism
of the space $\C^\N$ of all complex sequences.
Given $z\in \C$, we let
$V^z$ be the generalized
eigenspace of~$S$,
consisting of all $f\in \C^\N$ such that
$(S-z)^kf=0$ for some $k\in \N$.
Then $h_{k,z}\in V^z$
for each~$z\in \C^\times$ 
and $k\in \N_0$,
by a trivial induction:
If $k=0$,
we have
$((S-z)h_{0,z})(n)=z^{n+1}-z^{n+1}=0$.
If $k\geq 1$, then
$((S-z)h_{k,z})(n)=(n+1)^kz^{n+1}-n^kz^{n+1}
=((n+1)^k-n^k)z^{n+1}$.
Hence
$(S-z)h_{k,z}$
is a linear combination of
$h_{0,z},\ldots, h_{k-1,z}$,
each of which is in~$V^z$ by induction.\\[2.5mm]
Step~4: The sum $\sum_{z\in \C^\times}V^z$
of generalized eigenspaces
being direct, the\linebreak
assertion follows from Steps~1--3.
\end{proof}
Also the following simple lemma
is useful for the proof
of Theorem~\ref{mainres}
and later arguments.
It is a variant of~\cite[Lemma~10.1]{BGN}.
%
%
\begin{la}\label{basictool}
Let $E$ be a locally convex space over
$\K\in \{\R,\C\}$
and $E_0\sub E$ be a vector subspace
equipped with a locally convex vector topology making the
inclusion map $\iota\colon E_0\to E$ continuous.
Let $k\in \N_0$
and $f\colon U\to E$ be a map
on an open set $U\sub \K$, such that $f(U)\sub E_0$.
Then the following holds:
\begin{itemize}
\item[\rm(a)]
If the corestriction $f|^{E_0}\colon U\to E_0$ is $C^k_\K$, then
$f$ is $C^k_\K$ and $f^{(j)}(U)\sub E_0$ for all $j\in \N$ such that
$j\leq k$.
\item[\rm(b)]
If $E_0$ carries the topology induced by~$E$,
then $f|^{E_0}\colon U\to E_0$ is $C^k_\K$
if and only if~$f$
is $C^k_\K$ and $f^{(j)}(U)\sub E_0$ for all $j\in \N$ such that
$j\leq k$.
\end{itemize}
\end{la}
\begin{proof}
(a) By the Chain Rule, $f=\iota\circ f|^{E_0}$
is $C^k_\K$ with $f^{(j)}(x)=\iota\circ f^{(j)}(x)$,
from which the assertion follows.

(b) In view of~(a), we only need to prove
sufficiency of the described condition.
We proceed by induction on~$k$.
The case~$k=0$ being trivial, assume that $k\geq 1$.
If $f$ satisfies the described condition, then
$\frac{f|^{E_0}(y)-f|^{E_0}(x)}{y-x}
=\frac{f(y)-f(x)}{y-x}\to f'(x)\in E_0$ as $y\to x$,
showing that $f|^{E_0}$ is differentiable with
$(f|^{E_0})'=f'|^{E_0}$.
Since
$f'|^{E_0}$ is $C^{k-1}_\K$
by induction,
we see that $f$ is~$C^k_\K$.
\end{proof}
{\bf Proof of Theorem~\ref{mainres}.}
Consider the map
$f \colon U \to \C^\N$,
$f(z):=(e^{zn})_{n\in \N}$.
The $n$-th component $U\to \C$, $z\mto e^{zn}$
of~$f$ being $C^\infty_\C$ for each~$n$,
also the map $f$ into the direct product~$\C^\N$
is $C^\infty_\C$ (see \cite[Chapter~1]{GaN}
or \cite[Lemma~10.2]{BGN}),
with $f^{(j)}(z)=(n^je^{nz})_{n\in \N}$
and thus $f^{(j)}(z)=h_{j,e^z}$.
If $j\leq k$, then $h_{j,e^z}\in E_k$
and thus $f^{(j)}(U)\sub E_k$.
Thus $f_k=f|^{E_k}$ is $C^k_\C$,
by Lemma~\ref{basictool}\,(b).
However, $f^{(k+1)}(z)=h_{k+1,e^z}\not\in E_k$
by Definition of~$E_k$ and
Proposition~\ref{linind}.
Hence $f_k=f|^{E_k}$ is not $C^{k+1}_\C$,
by Lemma~\ref{basictool}\,(b).
If~$k=0$, then
$f_0$ is not $C^1_\C$ by the preceding.
However,
$\int_{\partial \Delta}f_0(\zeta)\,d\zeta=
\int_{\partial \Delta} f(\zeta)\,d\zeta=0$
for each triangle $\Delta\sub U$
(by ``(b)$\impl$(e)''
in the introduction),
because $f$ is~$C^\infty_\C$.
Hence~$f_0$ satisfies
property~(e) from the introduction
and hence~$f_0$
is weakly analytic.\,\Punkt
%
%
%
%
%
%
%
%
%
%
%
%
\section{Example of a complex analytic map
which \hspace*{.8mm}is not locally bounded}\label{secnbd}
Let $E:=\ell^1(\N,\C)$ be the space of absolutely summable
complex sequences with its usual norm $\|.\|_1$.
We define
\[
g\colon E\to\C\,,\quad g(x):=\sum_{k=1}^\infty 2^k\, (x_k)^{2k}
\]
for $x=(x_k)_{k\in \N}\in E$ and $f\colon E\to\C^\N$,
$f(x):=(g(nx))_{n\in \N}$.
We now show,
using the notion of complex analyticity
described in~{\rm (ii)}
in the introduction:
%
%
\begin{prop}\label{propunb}
The map $f\colon E\to\C^\N$ is complex
analytic.
It is unbounded on each non-empty
open subset of~$E$.
\end{prop}
As a preliminary, we discuss~$g$.
%
\begin{la}\label{lemong}
$g$ is complex analytic.
It is unbounded
on $\wb{B}_2^E(x)$
for all $x\in E$.
\end{la}
\begin{proof}
Since the partial sums
$g_n\colon E\to \C$, $x\mto \sum_{k=1}^n 2^k\, (x_k)^{2k}$
are polynomials in the point evaluations $x\mto x_k$
(which are continuous linear functionals)
and hence complex analytic,
$g$ will be complex analytic
if we can show that each $x\in E$ has
an open neighbourhood~$U$ such that $(g_n|_U)_{n\in \N}$
converges
uniformly (see \cite[Proposition~6.5]{BaS}).
There is $m\in \N$ such that $|x_k|<\frac{1}{4}$ for all
$k\geq m$. Set $U:=B_{1/4}^E(x)$.
Then $|y_k|<\frac{1}{2}$ for all $y\in U$
and thus $\sum_{k=m}^\infty\sup_{y\in U}|2^k\,(y_k)^{2k}|
\leq \sum_{k=m}^\infty 2^{-k}<\infty$,
which entails uniform convergence on~$U$.\\[2.5mm]
Given $x\in E$ and $N\in \N$, there exists
$m\in \N$ such that $|2x_m|<1$
and $2^m\geq N+|g(x)|+1$.
Set $y:=(y_k)_{k\in \N}$,
where $y_k:=x_k$ if $k\not=m$,
and $y_m:=1$. Then $\|y-x\|_1\leq 2$
and $|g(y)|=|g(x)+2^m-2^mx_m^{2m}|\geq 2^m-|g(x)|-2^m|x_m|^{2m}
\geq N$.
Thus $g$ is unbounded on~$\wb{B}_2^E(x)$.
\end{proof}
{\bf Proof of Proposition~\ref{propunb}.}
Given $x\in E$ and a neighbourhood $U$ of~$x$ in~$E$,
there exists $n\in \N$ such that
$\wb{B}_{2/n}^E(x)\sub U$.
Then $\wb{B}^E_2(nx)=n\wb{B}_{2/n}^E(x)\sub n U$,
whence $g(nU)\sub \C$ is unbounded,
by Lemma~\ref{lemong}.
Therefore the projection of $f(U)$
on the $n$-th component is unbounded
and thus $f(U)$ is unbounded.\,\Punkt
%
%
%
%
%
%
%
\section{Mappings to products and real analyticity}\label{secprod}
We describe a map $f=(f_n)_{n\in \N}\colon \R\to\R^\N$
which is not real analytic although each of its components
$f_n\colon \R\to\R$ is real analytic
(see also \cite[Chapter~II, Example~6.8]{Gro}
for a very similar example).
\begin{example}
For each $n\in\N$, the map
$f_n\colon \R\to\R$, $f_n(t):=\frac{1}{1+(nt)^2}$
is real analytic and its Taylor series
around~$0$ has radius of convergence~$\frac{1}{n}$.
It follows that the Taylor series around~$0$
of the smooth map
\[
f\colon \R\to\R^\N\,,\quad f(t):=(f_n(t))_{n\in \N}
\]
has radius of convergence~$0$,
and thus $f$ is not real analytic.
\end{example}
%
%
%
%
%
\section{Example of a map which is not real
analytic although it admits Taylor expansions}\label{secnoexp}
We describe a smooth map $f\colon \R\to E$ to a
suitable real locally convex space~$E$ which is
given locally by its Taylor series around
each point, but which does not admit a complex analytic
extension and hence fails to be real analytic.
We observe first that such a pathology cannot
occur if $E$ is Mackey complete.
Real analyticity is
understood as in~(iii) in the introduction.
%
%
\begin{la}\label{lanopatho}
Let $E$ be a Fr\'{e}chet space, $F$ be a Mackey complete
locally convex space and $f\colon U\to F$ be
a smooth map on an open subset $U\sub E$
which is locally given by its Taylor series
around each point $($with pointwise convergence$)$.
Then $f$ is real analytic.
\end{la}
\begin{proof}
The hypothesis means that $f$
is a real analytic map in the sense
of~\cite{BaS}. If~$F$ is sequentially complete,
\cite[Theorem~7.1]{BaS} provides
a complex analytic extension
for~$f$ into~$F_\C$
(because~$E$ is a Fr\'{e}chet space),
and thus~$f$ is real analytic in the desired sense.\\[2.5mm]
If~$F$ is merely Mackey complete,
we know from the preceding that $f$ is real analytic
as a map into the completion~$\wt{F}$ of~$F$.
Let $g\colon V\to \wt{F}_\C$ be a complex analytic
extension of~$f$ to an open neighbourhood
$V\sub E_\C$ of~$U$.
Given $x\in U$, let $W_x\sub E_\C$ be a balanced,
open $0$-neighbourhood with $x+W\sub V$.
Then $g(x+w)=\sum_{n=0}^\infty\frac{\delta^n_xf(w)}{n!}$
for each $w\in W_x$,
where $\delta^n_xf(w):=d^nf(x,w,\ldots, w)$.
Given $w\in W_x$,
there exists $t>1$ such that
$tx\in W_x$. Then
$\sum_{n=0}^\infty t^n \frac{\delta^n_xf(w)}{n!}
=\sum_{n=0}^\infty \frac{\delta^n_xf(tw)}{n!}$
converges in $\wt{F}_\C$,
whence $(\frac{\delta^n_xf(tw)}{n!})_{n\in \N_0}$
is a bounded sequence in~$F_\C$.
Since $\sum_{n=0}^\infty t^{-n}<\infty$,
the second characterization (M2)
of Mackey completeness in the introduction
shows that
$\sum_{n=0}^\infty \frac{\delta^n_xf(w)}{n!}
=\sum_{n=0}^\infty t^{-n} \frac{\delta^n_xf(tw)}{n!}$
converges in~$F_\C$.
Thus $g(x+w)\in F_\C$.
So,
after replacing $V$ by $\bigcup_{x\in U}(x+W_x)$,
we may assume that $g(V)\sub F_\C$.
Then $g\colon V\to F_\C$
is complex analytic
(see \cite[Proposition~1.5.18]{GaN})
and thus $f\colon U\to F$ is real analytic.
\end{proof}
Let $E$ be the space of
all sequences $x=(x_n)_{n\in \N}$
of real numbers
which have polynomial growth,
i.e.,
there exists $m\in \N$
such that the sequence
$(|x_n|n^{-m})_{n\in\N}$
is bounded.
We equip~$E$
with the
topology induced by~$\R^\N$.
%
%
\begin{prop}\label{propreal}
For~$E$ as before, the map
\[
f\colon \R\to E\,,\quad
f(t):=(\sin (nt))_{n\in \N}
\]
is smooth and $f(t)=\sum_{k=0}^\infty \frac{f^{(k)}(t_0)}{k!}(t-t_0)^k$
in~$E$, for all $t,t_0\in \R$.
However, $f$ is not real analytic.
\end{prop}
\begin{proof}
The map $g=(g_n)_{n\in \N}\colon \R\to\R^\N$,
$g(t):=f(t)$ is smooth,
since all components
$g_n\colon \R\to\R$, $g_n(t)=\sin(nt)$
are smooth.
We have
$g^{(2k)}(t)=(g_n^{(2k)}(t))_{n\in \N}
=(n^{2k}({-1})^k\sin (nt))_{n\in \N}\in E$
and
$g^{(2k+1)}(t)=(g_n^{(2k+1)}(t))_{n\in \N}
=(n^{2k+1}({-1})^k\cos(nt))_{n\in \N}\in E$
for each $k\in \N_0$ and $t\in \R$,
whence $f=g|^E$ is smooth
(by Lemma~\ref{basictool}\,(b)).\\[2.5mm]
Given $t_0\in\R$, we have
$g_n(t)
=\sum_{k=0}^\infty\frac{g_n^{(k)}(t_0)}{k!}(t-t_0)^k$
for each $n\in \N$ and $t\in \R$.
Since~$E$ is equipped with the topology induced by
$\R^\N$, it follows that
$f(t)=\sum_{k=0}^\infty\frac{f^{(k)}(t_0)}{k!}(t-t_0)^k$
for all $t,t_0\in\R$ and thus~$f$
is given globally by its Taylor series around
each point.
The map $g$ is real
analytic, because
\[
h\colon \C\to \C^\N\,,\quad h(z)=(\sin(nz))_{n\in \N}
\]
is a complex analytic extension of~$g$.
If~$f$ was complex analytic,
then we would have $h(V)\sub E_\C$
for some open neighbourhood $V$ of~$\R$
in~$\C$.
Then $it\in V$ for $t>0$
sufficiently small. However
\[
|\sin(int)| \; =\; |\sinh(nt)|\; \geq \; \frac{1}{4}e^{nt}
\]
for large~$n$ and hence
the sequence $h(it)=(\sin(int))_{n\in \N}$
does not have polynomial growth.
Thus $h(it)\not\in E_\C$, a contradiction.
\end{proof}
%
%
%
%
%
%
%
\section{Consequences for free locally
convex spaces}\label{secfreec}
Given $\K\in \{\R,\C\}$
and a completely regular topological
space~$M$, there exists
a Hausdorff locally convex topological
$\K$-vector
space $L(M,\K)$ and a continuous
map $\eta\colon M\to L(M,\K)$
with the following properties
(see \cite{Mar}):
\begin{itemize}
\item[(a)]
Algebraically,
$(L(M,\K),\eta)$ is the free $\K$-vector
space ($\isom \K^{(M)}$)
over the set~$M$;
\item[(b)]
For each continuous map
$\alpha \colon M\to E$ to a locally
convex topological $\K$-vector space~$E$,
there exists a unique continuous $\K$-linear map\linebreak
$\wb{\alpha}\colon L(M,\K)\to E$
such that $\wb{\alpha}\circ\eta=\alpha$.
\end{itemize}
$(L(M,\K),\eta)$ is determined by these
properties up to canonical isomorphism;
it is called the \emph{free locally convex
topological $\K$-vector space over~$M$.}
%
\begin{numba}\label{reducecx}
It is easy to see that $L(M,\C)$
is the complexification of~$L(M,\R)$
(by checking the universal property
for $L(M,\R)_\C$).
Hence $L(M,\R)$ is complete
(resp., sequentially complete,
resp., Mackey complete)
if and only if so is $L(M,\C)=L(M,\R)_\C$.
\end{numba}
%
%
\begin{prop}\label{consfree}
If $M\sub \C$ is a subset with non-empty interior~$M^0$,
then neither $L(M,\R)$ nor
$L(M,\C)$ is Mackey complete.
\end{prop}
\begin{proof}
Since~$\C$ is homeomorphic to the disk
$B^\C_1(1)$, after replacing~$M$
with a homeomorphic copy we may assume
that $0\not\in M$.
By \S\ref{reducecx},
we only need to show that $L(M,\C)$
is not Mackey complete.
Let $U\sub \C$ be a non-empty open set
with compact closure $\wb{U}$,
such that $e^{\wb{U}}\sub M$.
Set $E\!:=\!\Spann_\C\{h_{0,z}\colon z\in M\}$,
with $h_{0,z}$
as in~(\ref{defhs}).
Write~$E_0$ for~$E$, equipped with
the topology induced by the
direct product~$\C^\N$.
Let $\cO$ be the finest locally convex topology
on~$E$ such that $\eta\colon M\to E$,
$z\mto h_{0,z}$ is continuous.
Since the topology on~$E_0$ is Hausdorff
and makes~$\eta$ continuous, it
follows that
$\iota\colon (E,\cO)\to E_0$, $x\mto x$
is continuous and
$\cO$ is Hausdorff.
Since $(h_{0,z})_{z\in M}$ is a basis
for~$E$ by Proposition~\ref{linind},
it follows that $(E,\cO)$
is the free complex locally convex space
$L(M,\C)$
over~$M$
(together with~$\eta$).
Now consider
\[
g\colon \wb{U}\to E\,,\quad g(z)\,:=\,h_{0,e^z}
\,=\, (e^{nz})_{n\in \N}\,.
\]
Then $g=\eta\circ \exp|_{\wb{U}}$
is continuous, entailing
that $g(\wb{U})$ is compact and hence
bounded in $(E,\cO)$.
The restrictions $\lambda_n:=\pi_n|_{E}\to \C$
of the projections $\pi_n\colon \C^\N\to\C$,
$(x_k)_{k\in \N}\mto x_n$
are continuous linear on~$(E,\cO)$
and separate points.
Furthermore, $\lambda_n\circ g|_U\colon U\to\C$,
$z\mto e^{nz}$ is complex
analytic for each $n\in \N$.
Hence, if $(E,\cO)$ was Mackey complete,
then $g|_U\colon U\to E$ would be complex
analytic (by \cite[Theorem~1]{Gr2}).
But then also the map $f_0=\iota\circ g|_U\colon U\to E_0$
considered in Theorem~\ref{mainres}
would be complex analytic,
which it is not: contradiction.
Hence $L(M,\C)$ is not Mackey complete.
\end{proof}
%
%
\begin{rem}\label{remusp}
In the literature,
one finds various results concerning
$L(M,\R)$ and its completion,
which can be realized as a certain space
of measures (see \cite{Usp}, also~\cite{Flo}).
A result from~\cite{Usp} is of particular relevance:\\[2.5mm]
\emph{$L(M,\R)$ is complete if and only
if~$M$ is Dieudonn\'{e} complete\footnote{That is,
the largest uniformity on~$M$
compatible with its topology
is complete.}
and its compact subsets are finite.}\\[2.5mm]
Hence
$L(M,\R)$ and $L(M,\C)$
are non-complete in
the situation of Proposition~\ref{consfree}.
Our proposition provides
the additional information that
$L(M,\R)$ and $L(M,\C)$
are not sequentially complete
either, nor Mackey complete.
\end{rem}
Let us close with some
observations concerning
the free (not necessarily locally convex\,!)
topological $\K$-vector space $V(M,\K)$
over a completely regular
topological space~$M$ (obtained
by replacing the topology on $L(M,\K)$
with the finest vector topology making~$\eta$
continuous).\\[3mm]
To start with, let $M\sub \C$
be a compact set with non-empty interior.
Let $\K\in \{\R,\C\}$.
Then $V(M,\K)$ is complete,
by~\cite{AaK}.
If $V(M,\K)$ was locally convex,
we would have $L(M,\K)=V(M,\K)$
and so $L(M,\K)$ would be complete,
contrary to Proposition~\ref{consfree}.
We conclude: \emph{$V(M,\K)$
is not locally convex.}\\[2.5mm]
This argument can be generalized further.
To this end, recall that a Hausdorff
topological space~$M$ is said to be
a \emph{$k_\omega$-space}
if there exists a sequence $K_1\sub K_2\sub\cdots$
of compact subsets of~$M$ with union~$M$
(a ``$k_\omega$-sequence'')
such that $M=\dl\, K_n$\vspace{-.6mm}
as a topological space\,\footnote{Thus,
a set $A\sub M$ is closed if and only if
$A\cap K_n$ is closed in~$K_n$ for each
$n\in \N$.}
(see, e.g., \cite{GGH}
and the references therein for further information).
For instance, every $\sigma$-compact\linebreak
locally compact topological
space is a $k_\omega$-space.
Each $k_\omega$-space
is normal
(by \cite[Proposition~4.3\,(i)]{Han})
and hence
completely regular,
ensuring that both
$L(M,\K)$ and $V(M,\K)$ are defined.
We show:
\begin{prop}
If $M$ is a non-discrete
$k_\omega$-space,
then neither $V(M,\R)$ nor $V(M,\C)$
is locally convex.
\end{prop}
\begin{proof}
Let $K_1\sub K_2\sub \cdots$
be a $k_\omega$-sequence
for~$M$. If each $K_n$ was finite,
then $K_n$ would be discrete
and hence also $M=\dl\,K_n$\vspace{-.4mm}
would be discrete,
contradicting the hypothesis.
Therefore some~$K_n$ is infinite,
whence $L(M,\R)$
(and hence also $L(M,\C)$)
is not complete,
by Uspenski\u{\i}'s result recalled
in Remark~\ref{remusp}.
Since $V(M,\K)$
is complete
by the next lemma,
we see that it cannot
coincide with $L(M,\R)$
and hence cannot be locally convex.
\end{proof}
\begin{la}
Let $M$ be a $k_\omega$-space
and $\K\in \{\R,\C\}$.
Then
$V(M,\K)$
is a $k_\omega$-space
and hence complete.
\end{la}
\begin{proof}
Since each abelian topological
group which is a $k_\omega$-space is complete~\cite{Rai},
we only need to show
that $V(M,\K)$
is a $k_\omega$-space.
To this end,
it is convenient to identify~$M$ via~$\eta$
with a subset of $V(M,\K)$.
We pick a $k_\omega$-sequence
$(K_n)_{n\in \N}$
for~$M$.
Define
$L_n:=\wb{B}^\K_n(0)
\cdot (K_n+\cdots+K_n)$
(with $2^n$ summands) for $n\in \N$.
Then $L_1\sub L_2\sub\cdots$ is a sequence
of compact
subsets of~$V(M,\K)$,
with union~$V(M,\K)$.
The topology~$\cO$
making $V(M,\K)$
the direct limit topological space
$\dl\, L_n$\vspace{-.4mm}
is finer
than the original
topology and makes~$V(M,\K)$
a $k_\omega$-space.
We now write~$W$ for~$V(M,\K)$,
equipped with the topology~$\cO$.
Because the inclusion map $\iota\colon M\to W$
restricts to a continuous
map on~$K_n$ for each~$n$
(since we can pass over~$L_n$),
we see that $\iota$ is continuous
(as $M=\dl\,K_n$).\vspace{-.5mm}
To complete the proof,
it only remains to show
that~$\cO$ is a vector topology.
Since $W\times W=\dl\,(K_n\times K_n)$\vspace{-.5mm}
(see \cite[Proposition~3.3]{DIR}),
the addition map $\alpha\colon W\times W\to W$
will be continuous
if $\alpha|_{K_n\times K_n}$
is continuous
for each~$n$. Since $K_n+K_n\sub K_{n+1}$
and~$W$ induces the same topology on
$K_n$ and on~$K_{n+1}$ as $V(M,\K)$,
continuity of
$\alpha|_{K_n\times K_n}$
follows
from the continuity
of the addition map
$V(M,\K)\times V(M,\K)\to V(M,\K)$.
Likewise, since $\K\times W=\dl\, \wb{B}^\K_n(0)\times K_n$\vspace{-.4mm}
and $\wb{B}^\K_n(0)K_n\sub K_{n^2}$,
we deduce from the continuity
of the scalar multiplication map
$\K\times V(M,\K)\to V(M,\K)$
that also the scalar multiplication
$\K\times W\to W$
is continuous.
\end{proof}
Similar
arguments show
that the free
topological group and
the free abelian topological
group over a $k_\omega$-space
are $k_\omega$-spaces
(see~\cite{MMO}).
{\footnotesize
{\bf Helge Gl\"{o}ckner}, TU Darmstadt, Fachbereich Mathematik AG~5,
Schlossgartenstr.\,7,\\
64289 Darmstadt, Germany.
\,E-Mail:
\,{\tt gloeckner@mathematik.tu-darmstadt.de}}
\end{document}